\documentclass{amsproc}

\usepackage[T1]{fontenc}
\usepackage[francais]{babel}
\usepackage{pstricks}
\usepackage{graphicx,psfrag,pstcol,pst-grad, pst-eucl, pstricks-add}
\usepackage{amssymb}
\usepackage{amsmath}
\usepackage{amsthm}
\usepackage{amscd}
\usepackage{amsfonts}
\usepackage{stmaryrd}
\usepackage{pb-diagram}
\usepackage{epic,eepic,epsfig}
\usepackage{a4wide}
\usepackage{enumerate}
\usepackage{nextpage}
\usepackage{fancyhdr}
\usepackage{epsf}

\newtheorem{theorem}{Th\'eor\`eme}[section]
\newtheorem{lemma}[theorem]{Lemme}
\newtheorem{conjecture}[theorem]{Conjecture}
\newtheorem{corollary}[theorem]{Corollaire}
\newtheorem{proposition}[theorem]{Proposition}

\theoremstyle{definition}

\theoremstyle{remark}
\newtheorem{remark}[theorem]{Remarque}

\numberwithin{equation}{section}

\newcommand{\Jac} {\mathop{\mathrm{Jac}}}
\newcommand{\hFplus} {\mathop{h_{\mathrm{F}^+}}}
\newcommand{\codim} {\mathop{\mathrm{codim}}}
\newcommand{\Card} {\mathop{\mathrm{Card}}}

\begin{document}

\title{Bornes sur le nombre de points rationnels des courbes : en qu\^ete d'uniformit\'e}

\author{Fabien Pazuki}
\address{University of Copenhagen, Institute of Mathematics, Universitetsparken 5, 2100 Copenhagen, Denmark, and Universit\'e de Bordeaux, IMB, 351, cours de la Lib\'eration, 33400 Talence, France.}
\curraddr{University of Copenhagen, Institute of Mathematics, Universitetsparken 5, 2100 Copenhagen, Denmark}
\email{fpazuki@math.ku.dk}
\thanks{L'auteur remercie Ga\"el R\'emond pour ses commentaires pr\'ecis et pr\'ecieux. Merci \`a Amos Turchet pour ses remarques. L'auteur est soutenu par la Chaire DNRF Niels Bohr de Lars Hesselholt et les Projets ANR-14-CE25-0015 Gardio et ANR-17-CE40-0012 Flair.}

%    General info
\date{}

\dedicatory{}


\begin{abstract}
L'objectif de ce texte est de montrer en d\'etail comment un r\'esultat conjectural de minoration de hauteur de type Lang-Silverman permet d'obtenir une borne explicite sur le nombre de points rationnels d'une courbe de genre $g\geq 2$ sur un corps de nombres. La borne est uniforme en la hauteur de la courbe.
\end{abstract}

\maketitle

\section{Introduction}

Soit $C$ une courbe projective lisse de genre $g\geq 2$ d\'efinie sur un corps de nombres $K$. La conjecture de Mordell, qui est depuis 1983 un th\'eor\`eme c\'el\`ebre de Faltings \cite{L}, donne la finitude de l'ensemble des points $K$-rationnels $C(K)$. Une question naturelle a alors \'emerg\'e : \`a quel point le cardinal $\#C(K)$ d\'epend-il de la courbe $C$ ? Est-il possible de produire un majorant explicite de $\#C(K)$ et \`a quel point ce majorant peut-il \^etre uniforme ? C'est \`a cette question d'uniformit\'e du majorant que ce texte est consacr\'e.

On sait produire depuis R\'emond \cite{W} des bornes totalement explicites et inconditionnelles. Citons une forme plus r\'ecente de son r\'esultat, donn\'ee dans \cite{X}.

\begin{theorem}
Soit $C$ une courbe lisse, projective, g\'eometriquement connexe et de genre $g\geq 2$ d\'efinie sur un corps de nombres $K$. On fixe un plongement projectif de la vari\'et\'e jacobienne $Jac(C)$ relatif \`a la puissance seizi\`eme du fibr\'e associ\'e \`a un translat\'e sym\'etrique du diviseur th\^eta. Alors on a $$\#C(K)\leq (2^{38+2g}\cdot[K:\mathbb{Q}]\cdot g \cdot \max\{1, h_{\theta}\})^{(m_K+1)\cdot g^{20}},$$ o\`u $h_{\theta}$ est la hauteur th\^eta de $Jac(C)$ dans ce plongement et $m_K$ est le rang de Mordell-Weil de Jac(C)(K).
\end{theorem}

On cherche \`a savoir dans quelle mesure on peut se passer de la d\'ependance en la jacobienne de la courbe dans l'expression du majorant. L'uniformit\'e la plus forte est conjectur\'ee par Caporaso, Harris, Mazur en page 2 de \cite{E} :

\begin{conjecture}\label{uniforme 1}
Soit $g\geq 2$ un entier naturel et soit $K$ un corps de nombres.
Il existe une quantit\'e $c_0(g,K)>0$ telle que pour toute courbe projective lisse $C$ de genre $g$ d\'efinie sur $K$, on a $$\#C(K)\leq c_0(g,K).$$
\end{conjecture} 

Ce dernier \'enonc\'e est impliqu\'e par une conjecture tr\`es g\'en\'erale de Lang qui s'\'enonce comme suit : les points rationnels sur une vari\'et\'e de type g\'en\'eral ne sont pas Zariski denses. La d\'eduction se fait via un proc\'ed\'e dit de corr\'elation qui est bien expliqu\'e dans \cite{E}. On consultera aussi avec profit les textes \cite{A, S} sur le sujet. Ce niveau d'uniformit\'e maximal peut sembler difficile \`a atteindre. On pourra essayer dans un premier temps\footnote{Ajout\'e apr\`es arbitrage : une preuve de la conjecture \ref{uniforme 2} est annonc\'ee par Dimitrov, Gao, Habegger dans \cite{K} !} de r\'epondre \`a la question classique suivante (trouv\'ee par exemple dans \cite{H} en introduction).

\begin{conjecture}\label{uniforme 2}
Soit $g\geq 2$ un entier naturel et soit $K$ un corps de nombres.
Il existe une quantit\'e $c_1(g,K)>0$ telle que pour toute courbe projective lisse $C$ de genre $g$ d\'efinie sur $K$, on a $$\#C(K)\leq c_1(g,K)^{m_K+1},$$ o\`u $m_K$ est le rang de Mordell-Weil de la jacobienne de $C$ sur $K$.
\end{conjecture} 

Citons aussi une version plus souple de la conjecture \ref{uniforme 2} et trouv\'ee sous forme de question dans \cite{R} page 223 (voir aussi le premier paragraphe de la page 234 de \cite{Q}):

\begin{conjecture}\label{uniforme 3}
Soit $g\geq 2$ un entier naturel, soit $K$ un corps de nombres, soit $m$ un entier naturel.
Il existe une quantit\'e $c_2(g,K,m)>0$ telle que pour toute courbe projective lisse $C$ de genre $g$ d\'efinie sur $K$, dont le rang de Mordell-Weil de la jacobienne de $C$ sur $K$ est $m$, on a $$\#C(K)\leq c_2(g,K,m).$$ 
\end{conjecture} 

David, Nakamaye et Philippon donnent dans le th\'eor\`eme 3.6 du texte \cite{H} des familles infinies de courbes d\'efinies sur $\mathbb{Q}$ et v\'erifiant la conjecture \ref{uniforme 2} (des familles infinies de courbes sont propos\'ees pour chaque $g\geq 4$ fix\'e). Notons qu'il est classique de penser que la d\'ependance en le corps $K$ des quantit\'es  $c_{0}, c_1, c_2$ ne fait intervenir que le degr\'e de $K$ sur $\mathbb{Q}$.

De nouveaux progr\`es, suivant une approche de ces questions bas\'ee sur l'int\'egration $p$-adique, concernent les courbes $C$ telles que la jacobienne $\Jac(C)$ poss\`ede un rang de Mordell-Weil sur $K$ petit par rapport au genre de $C$. Initi\'ee par des id\'ees de Chabauty, puis de Coleman, la strat\'egie $p$-adique de majoration du cardinal de $C(K)$ fait voir le jour \`a un r\'esultat d'uniformit\'e par Stoll \cite{Z} qui traite des courbes hyperelliptiques. Cette attaque a \'et\'e g\'en\'eralis\'ee \`a toutes les courbes (pour $g\geq 3$) dans \cite{N} par Katz, Rabinoff, Zureick-Brown, dont voici le r\'esultat :

\begin{theorem}\label{KaRaZu}
Soit $d\geq 1$ et soit $g\geq 3$ des entiers. Il existe une constante $c_3(g,d)>0$ telle que pour tout corps de nombres $K$ de degr\'e $d$ et pour toute courbe projective lisse $C$ de genre $g$, d\'efinie sur $K$ et telle que le rang $m_K$ de la jacobienne $\Jac(C)$ sur $K$ v\'erifie $m_K\leq g-3$, on a $\#C(K)\leq c_3(g,d)$.
\end{theorem}

Notons que la quantit\'e $c_3(g,d)$ peut \^etre explicit\'ee, comme cela est propos\'e dans les textes \cite{Z} et \cite{N}. Lorsqu'elle est applicable, la m\'ethode de Chabauty-Coleman fournit une d\'ependance polynomiale en le param\`etre $g$, voire quasi-lin\'eaire.

\begin{remark}
Le th\'eor\`eme \ref{KaRaZu} pourra donc \^etre vu comme un pas vers la conjecture \ref{uniforme 1}, ou comme un pas vers les conjectures \ref{uniforme 2} ou \ref{uniforme 3}. En effet la conjecture \ref{uniforme 2} (ou la conjecture \ref{uniforme 3}) combin\'ee \`a une in\'egalit\'e du type $m_K \leq c_4(g)$ pour une constante $c_4(g)>0$ ne d\'ependant que de $g$ impliquent la conjecture \ref{uniforme 1} \textit{restreinte aux courbes dont la jacobienne v\'erifie $m_K\leq c_4(g)$} bien entendu.
\end{remark}

\begin{remark}
La conjecture \ref{uniforme 2} ou la conjecture \ref{uniforme 3} coupl\'ee \`a une borne uniforme sur le rang de la forme $m_K\leq c(g,K)$ implique bien entendu la conjecture \ref{uniforme 1}. La question de savoir si le rang des vari\'et\'es ab\'eliennes de dimension fix\'ee $g\geq1$ sur un corps de nombre fix\'e $K$ est born\'e uniform\'ement est un probl\`eme encore largement ouvert. De r\'ecents travaux \cite{C} tendent \`a sugg\'erer qu'une borne uniforme est plausible pour les courbes elliptiques sur $\mathbb{Q}$.
\end{remark}

Notre modeste contribution ici est de montrer en d\'etail qu'une conjecture de minoration de hauteur dans l'esprit Lang-Silverman implique aussi la conjecture \ref{uniforme 2}, sans besoin d'hypoth\`ese sur le rang. C'est une id\'ee naturelle apr\`es la lecture de \cite{J}, qui traite cependant uniquement le cas des jacobiennes isog\`enes \`a des produits de courbes elliptiques (voir l'hypoth\`ese $(*)$ page 110 de \cite{J}). On traite le cas g\'en\'eral ici en suivant le chemin euclidien trac\'e par R\'emond dans ses articles \cite{W, X}, et en ajoutant comme ingr\'edient une minoration \`a la Lang-Silverman. La d\'eduction de la conjecture \ref{uniforme 2} s'obtient alors moyennant quelques efforts techniques suppl\'ementaires. Ce faisant, nous am\'eliorons donc la proposition 1.9 page 23 de \cite{T} qui est bas\'ee sur une version obsol\`ete de Lang-Silverman. Voici une version plus moderne de cette conjecture.

\begin{conjecture}(\cite{Y})\label{LangSilV3}
Pour tout couple d'entiers naturels $(d,g)$ il existe un r\'eel $c_{5}=c_{5}(d,g)>0$ tel que pour tout corps de nombres $K$ de degr\'e $d$, pour toute vari\'et\'e ab\'elienne $A/K$ de dimension $g$ et pour toute polarisation $L$ sur $A$, il existe une sous-vari\'et\'e ab\'elienne stricte $B$ de $A$ de degr\'e $\deg_L (B)\leq c_{5}\deg_L (A)$ et telle que pour tout point $P\in{A(K)}$, au moins l'une des deux assertions suivantes est v\'erifi\'ee :
\begin{enumerate}
\item il existe un entier $1\leq N\leq c_{5}$ tel que $[N]P\in B$,

\item
\[
\widehat{h}_{A,L}(P) \geq c_{5}^{-1}\, \max\Big\{\hFplus(A/K),\,1\Big\}\;,
\]
o\`u $\widehat{h}_{A,L}(.)$ est la hauteur de {N\'eron-Tate} associ\'ee \`a $L$ et $\hFplus(A/K)$ est la hauteur de Faltings de la vari\'et\'e $A/K$.
\end{enumerate}
\end{conjecture}

La formulation originale de cette conjecture pr\'esente uniquement la minoration du second cas, la conjecture \ref{LangSilV3} donne une description plus compl\`ete. Notons que la sous-vari\'et\'e $B$ ne d\'epend pas de $P$. L'utilisation de cette minoration en famille de la hauteur canonique sur les vari\'et\'es jacobiennes fournit alors le r\'esultat suivant.

\begin{theorem}\label{thm final}
Supposons vraie la conjecture $\ref{LangSilV3}$. Soit $K$ un corps de nombres de degr\'e $d$ et soit $g\geq 2$ un entier. Il existe une constante $c_{6}=c_{6}(g,K)>0$ telle que pour toute courbe $C$ de genre $g$ d\'efinie sur $K$, 
$$\#C(K)\leq c_{6}^{m_K+1},$$
o\`u $m_K$ d\'esigne le rang de Mordell-Weil de la jacobienne de $C$ sur $K$. On peut choisir   
$$c_6=c_{5}^{2g+2} 2^{120g} g^{14+g},$$ o\`u $c_5$ est la quantit\'e introduite dans la conjecture \ref{LangSilV3}.
\end{theorem}

On trouve des r\'esultats partiels en direction de la conjecture \ref{LangSilV3} dans les travaux \cite{U} en dimension 2 et \cite{F, P} en dimension g\'en\'erale. Ceci fournit donc le corollaire imm\'ediat suivant, non couvert par \cite{H, J}.

\begin{corollary}
 Les familles de courbes de genre $2$ satisfaisant aux hypoth\`eses du corollaire 1.14 de \cite{U} et de genre $g\geq 2$ d\'ecrites dans \cite{F, P} v\'erifient la conjecture \ref{LangSilV3} et sont donc des exemples o\`u la borne du th\'eor\`eme \ref{thm final} est d\'emontr\'ee inconditionnellement.
\end{corollary}

La section \ref{euclide} contient les arguments de comptage euclidien et la preuve du th\'eor\`eme \ref{thm final}, la section \ref{varianteLS} pr\'esente une variante de Lang-Silverman et les modifications que cela entra\^ine dans le comptage, la section \ref{DavPhi} est un appendice \'ecrit par Sinnou David et Patrice Philippon qui corrigent un lemme de leur article \cite{I} et fournissent par la m\^eme occasion un calcul utile \`a la comparaison de hauteurs donn\'ee en (\ref{comparaison}).

\section{Comptage euclidien}\label{euclide}

On donne dans cette section une preuve d\'etaill\'ee du fait qu'une version forte de Lang-Silverman implique une borne uniforme sur le nombre de points rationnels des courbes de genre $g\geq 2$ sur les corps de nombres. 

Dans toute la suite, on se place dans le cadre suivant. Soit $C$ une courbe projective lisse et g\'eom\'etriquement connexe, de genre $g\geq 2$ sur un corps de nombres $K$. On suppose tout d'abord que la courbe poss\`ede un point rationnel $P_0\in C(K)$, dans le cas contraire toute borne sup\'erieure positive donn\'ee sur le cardinal de $C(K)$ sera trivialement valide. La courbe $C$ est alors plong\'ee dans sa jacobienne $A=\Jac(C)$ par le plongement $j: P\to (P)-(P_0)$, et on d\'efinit un diviseur th\^eta par $\Theta=j(C)+...+j(C)$ o\`u on somme $g-1$ termes $j(C)$, diviseur fournissant une polarisation principale sur $\Jac(C)$. On consid\`ere la jacobienne plong\'ee dans $\mathbb{P}^{4^{2g}-1}$ via le fibr\'e $L$ correspondant \`a la puissance 16-i\`eme du fibr\'e associ\'e au translat\'e sym\'etrique de ce diviseur th\^eta. Par le th\'eor\`eme 5.8 page 118 de \cite{O}, on a $\deg_{\Theta}(C)=g$. On aura donc $\deg_L(A)=2^{4\dim(A)}\deg_{\Theta}(A)=16^g g!$ et $\deg_L(C)=2^{4\dim(C)}\deg_{\Theta}(C)=16g$. 

Si on note $\phi:A\to \mathbb{P}^{4^{2g}-1}$ le plongement associ\'e, la hauteur th\^eta de $A$ est alors d\'efinie par $h_\Theta(A)=h(\phi(0_A))$, o\`u $h(.)$ est la hauteur du vecteur coordonn\'ees comme donn\'ee dans \cite{X} page 762. On utilise aussi la hauteur de Faltings $\hFplus(A)$ d'une vari\'et\'e ab\'elienne $A$, normalis\'ee comme dans \cite{V}. On utilise enfin la classique hauteur de N\'eron-Tate d'un point alg\'ebrique sur une vari\'et\'e ab\'elienne $A$ polaris\'ee par un fibr\'e ample et symm\'etrique $L$, qui est une forme quadratique not\'ee $\widehat{h}_{A,L}$ ou simplement $\widehat{h}$ quand le contexte est clair.

Nous allons faire usage de la proposition 3.7 page 527 de \cite{W} dans une version valable uniquement pour les courbes et reprenant les r\'esultats du paragraphe baptis\'e ``Cas des courbes'' page 529 de \cite{W}, lequel explique comment am\'eliorer les valeurs des constantes dans le cas particulier qui nous int\'eresse ici. Cet \'enonc\'e est donc num\'eriquement plus efficace qu'une banale sp\'ecialisation de la proposition 3.7 page 527 au cas des courbes. 

\begin{proposition}(R\'emond)\label{remond1}
Soit $C$ une courbe de genre $g\geq 2$ plong\'ee dans sa jacobienne $A=\Jac(C)$ comme pr\'ec\'edemment dans un espace projectif ambiant $\mathbb{P}^n$ avec $n=4^{2g}-1$, soit $\Gamma$ un sous-groupe de $A(\overline{\mathbb{Q}})$ de rang fini $r$. Soit $c_{NT}$ une constante majorant la diff\'erence entre la hauteur de Weil et la hauteur de N\'eron-Tate sur $A$, relativement au plongement donn\'e. Soit $h_{1}$ la hauteur de la famille de polyn\^omes d\'efinissant l'addition sur $A$, relativement au plongement donn\'e. On d\'efinit alors $c_7$ par la formule 
$$c_{7}^{2}=2^{187+8g} g^{21}\max\{1, c_{NT}, h_1\},$$
et on a l'alternative suivante, o\`u au moins une des assertions est v\'erifi\'ee :

\begin{itemize}
\item ou 
$$\#\Big(C\cap \Gamma \Big)\leq (2^{170+8g} g^{12})^{r},$$

\item ou il existe un point $y\in{C\cap \Gamma}$ tel que $$\#\Big\{x\in{C\cap \Gamma}\;|\;\widehat{h}(x-y)\geq 4c_{7}^{2}\Big\}\;\,\leq \, (2^{170+8g} g^{12})^{r}.$$
\end{itemize}
\end{proposition}

\begin{proof}
Le c\oe ur de la preuve est bien entendu d\'etaill\'e dans \cite{W} pages 527-529. On demandera au lecteur de se r\'ef\'erer au texte original pour les arguments clefs de la preuve de cette proposition. Nous donnons ici les valeurs des quantit\'es correspondant \`a notre contexte. Tout d'abord comme $g\geq 2$, l'union des translat\'es de sous-vari\'et\'es ab\'eliennes contenues dans $C$ est vide, on a donc $Z_C=\emptyset$ dans la notation de \cite{W}. La quantit\'e $c_5$ apparaissant dans \cite{W} vient de la Proposition 3.6 page 526, mais est am\'elior\'ee par le paragraphe traitant le cas particulier des courbes page 529, on peut donc choisir comme valeur pour $c_7$ la quantit\'e donn\'ee par $c_{7}^{2}=2^{20}\Lambda^2 c_{NT} +2(n+1)\Lambda^{28/3} \max\{((n+1)D+1)(c_{NT}+3\log(n+1)), \Lambda^{-1} h_1\}$, o\`u $\Lambda = \max\{(\deg_L(C))^2, 2^{14} \deg_L(C)\}$ et $\deg_L(C)=16g$. En majorant l\'eg\`erement plus et en tenant compte de $n=4^{2g}-1$ et de $D=\max\{2^{14}, \deg_L(C)\}\leq 2^{14}g$, on voit que $2^{187+8g} g^{21}\max\{1, c_{NT}, h_1\}$ convient.

Le majorant figurant dans l'alternative est donn\'e par $(3(n+1)^2D^{12})^{r}$ dans \cite{W} page 529, ce qui fournit la valeur propos\'ee ici l\`a encore via $n=4^{2g}-1$ et $D=\max\{2^{14}, \deg_L(C)\}\leq 2^{14}g$.
\end{proof}

On donne \`a pr\'esent un lemme utile, une estimation sur la taille du rayon discriminant les points ``petits'' des points ``grands'' dans le comptage euclidien qui suivra.

\begin{lemma}\label{c5}
Soit $C$ une courbe de genre $g\geq 2$. Nous conservons le m\^eme cadre que pr\'ec\'edemment, un plongement th\^eta de la jacobienne $\Jac(C)$ dans $\mathbb{P}^{n}$, avec $n=4^{2g}-1$, correspondant \`a un translat\'e sym\'etrique du diviseur th\^eta (multipli\'e par $16$). Soit $h_{\Theta}(A)$ la hauteur th\^eta de la vari\'et\'e ab\'elienne $A=\Jac(C)$ dans ce plongement et consid\'erons la quantit\'e $c_{7}^2$ pour ce plongement de $A$ dans $\mathbb{P}^{n}$. Alors en utilisant les notations pr\'ec\'edentes 
$$c_{7}^{2}\leq 2^{196+10g}g^{22}\max\{1, h_{\Theta}(A)\}.$$
\end{lemma}

\begin{proof}
On part de la valeur explicite de la proposition pr\'ec\'edente
$$c_{7}^{2}=2^{187+8g} g^{21}\max\{1, c_{NT}, h_1\},$$
et on trouve dans \cite{I} les majorations (page 652, page 662 et page 665), 
$$h_{1}\leq 2(4^{g}-1)h_{\Theta}(A)\;,\;\;\;\;\;\;\;\; c_{NT}\leq (4^{g+2}h_{\Theta}(A)+6g\log2)\deg_L(C).$$
Comme plus haut, $C$ est plong\'ee dans sa jacobienne par un plongement th\^eta, on a alors dans $\mathbb{P}^{4^{2g}-1}$ la valeur $\deg_L(C)=16g$, o\`u $L$ correspond \`a la puissance seizi\`eme du fibr\'e associ\'e au diviseur th\^eta. Ceci implique
$$c_{7}^{2}\leq 2^{187+8g} g^{21}\max\{1, 2^{2g+9}g, 2^{2g+1}\} \max\{1,h_{\Theta}(A)\},$$
et on obtient donc
$$c_{7}^{2}\leq 2^{196+10g}g^{22}\max\{1, h_{\Theta}(A)\}.$$
\end{proof}

Nous sommes \`a pr\'esent en mesure de d\'emontrer le th\'eor\`eme \ref{thm final} de l'introduction.

\begin{proof} (du th\'eor\`eme \ref{thm final}) 
Soit $A$ la jacobienne de $C$ et soit $m_{K}$ le rang de $A(K)$.
On consid\`ere $C(K)$ comme un sous-ensemble de $A(K)$. Notons comme toujours $\widehat{h}$ la hauteur de N\'eron-Tate sur $A$ associ\'ee au plongement th\^eta dans $\mathbb{P}^{4^{2g}-1}$. Nous allons en fait borner $\#C(K)$ en majorant $\#(C\cap A(K))$. 

D'apr\`es la proposition \ref{remond1}, il suffit de traiter le cas o\`u il existe un point $y\in{C(K)}$ tel que $$\#\Big\{x\in{C(K)}\;|\;\widehat{h}(x-y)\geq 4c_{7}^{2}\Big\}\;\,\leq \,N,$$ o\`u $N$ est en fait un entier explicite. On divise les points rationnels en deux clans :

\begin{center}
\begin{tabular}{ll}
$S_1=$ & $\Big\{x\in{C(\overline{\mathbb{Q}})\cap A(K)\, \Big\vert\,\widehat{h}(x-y)\leq 4c_{7}^2}\Big\}$,
\\

$S_2=$ & $\Big\{x\in{C(\overline{\mathbb{Q}})\cap A(K)\,\Big\vert\, \widehat{h}(x-y)\geq 4c_{7}^2}\Big\}.$
\\
\end{tabular}
\end{center}

Notre entreprise de majoration de $\#C(K)$ commence par la remarque triviale $C(K)\subset S_1\cup S_2.$
On majore alors le cardinal des deux ensembles auxiliaires $S_1$ et $S_2$.

Invoquons tout d'abord la proposition \ref{remond1}, qui fournit directement 
\begin{equation}\label{S2}
\#S_2\leq (2^{170+8g} g^{12})^{m_K}.
\end{equation}

Pour borner le cardinal de $S_1$, il faut travailler un peu plus. On adapte ici deux techniques respectivement expos\'ees dans \cite{J} et dans \cite{W}  \`a notre situation. L'ensemble $(A(K)\otimes\mathbb{R}, \widehat{h})$ est un espace euclidien, on gardera les m\^emes notations pour l'image des points de $A(K)$ dans cet espace. Le fait que $\widehat{h}$ est une forme quadratique d\'efinie positive sur cet espace euclidien est classique (mais non trivial de prime abord, comme remarqu\'e par Cassels), une preuve est donn\'ee par exemple par la proposition B.5.3 de \cite{M} page 201.

Soit $B_o(P,\rho)$ la boule ouverte euclidienne de centre $P$ et rayon $\rho$. Une version du lemme g\'eom\'etrique \'el\'ementaire suivant est d\'emontr\'ee par R\'emond et figure au lemme 6.1 page 541 de \cite{W}.

\begin{lemma} \label{eucl}
Soit $r$ un entier et $\rho>0$ un r\'eel. Soit $S$ un sous-ensemble de $\mathbb{R}^r$ contenu dans une boule euclidienne de rayon $\rho$. Pour tout nombre r\'eel $\gamma\geq1$, on peut trouver $\lfloor (2\gamma+1)^r \rfloor$ boules de rayon $\rho/\gamma$, centr\'ees en des points de $S$ et telles que leur union recouvre tout $S$.
\end{lemma}

Nous pouvons illustrer la situation sur un dessin (en deux dimensions, certes). Les points noirs d\'esignent les points de $S$. On recouvre tout $S$ en choisissant des points centr\'es sur les ``paquets'' de points \'eventuellement regroup\'es et en faisant grandir le rayon des petits cercles jusqu'\`a $\rho/ \gamma$.

\begin{pspicture}(-4,-4)(5,6)

\pscircle(2,1){3.5}

\rput(2.4,1){$O$}
\psdot[dotstyle=+](2,1)

\psdots(2,3)(2,3.5)(2.2,3.6)(3.1,1,9)(3.6,3.3)(-1,0)(3,2.8)(1.2,1.3)(0.9,1.4)(4.2,-1.4)(4.2,-0.3)
\psdots(-0.5,0)(1,1)(2.2,3.4)(3.1,3.2)(3.6,0.8)(4.3,-1)%les centres
\pscircle(-0.5,0){0.8}
\pscircle(3.6,0.8){0.8}
\pscircle(1,1){0.8}
\pscircle(2.2,3.4){0.8}
\pscircle(3.1,3.2){0.8}
\pscircle(4.3,-1){0.8}

\psline{<->}(2,1)(2.4,-2.45)
\rput(2, -0.8){$\rho$}

\psline{<->}(3.6, 0.8)(4.1,0.18)
\rput(3.35, 0.45){$\rho/ \gamma$}

\rput(8,1.6){Illustration du}
\rput(8,1){lemme \ref{eucl}}

\end{pspicture}

On va appliquer ce lemme \`a $S=S_1\subset B_o(y,2c_{7})$. Ici $\rho=2c_{7}$ et $r=m_K$, donc en prenant $\gamma=\max\{1, 2c_{7}/\sqrt{c_{5}^{-1}\max\{1,\hFplus(A/K)\}}\}$ on obtient $$\frac{\rho}{\gamma}=\frac{2c_7}{\max\{1, 2c_{7}/\sqrt{c_{5}^{-1}\max\{1,\hFplus(A/K)\}}\}}=\min\{2c_7, \sqrt{c_{5}^{-1}\max\{1,\hFplus(A/K)\}} \}$$ ainsi que l'existence d'un ensemble fini $I\subset S_1$ de cardinal inf\'erieur \`a 
\begin{equation}\label{cardinal}
\lfloor (2\gamma+1)^r \rfloor\leq \left(1+\max\{2, 4c_{7}/\sqrt{c_{5}^{-1}\max\{1,\hFplus(A/K)\}}\}\right)^{m_K}
\end{equation}
tel que
$$S_1\subset\bigcup_{P_0\in I} B_o(P_0, \min\{2c_7, \sqrt{c_{5}^{-1}\max\{1,\hFplus(A/K)\}}\}).$$ Comme $S_1\subset C(K)$, on obtient  
$$S_1\subset\bigcup_{P_0\in I}B_o(P_0, \min\{2c_7, \sqrt{c_{5}^{-1}\max\{1,\hFplus(A/K)\}}\}) \cap C(K).$$
Or tout point $P$ de $B_o(P_0, \min\{2c_7, \sqrt{c_{5}^{-1}\max\{1,\hFplus(A/K)\}}\})$ v\'erifie en particulier l'in\'egalit\'e $$\hat{h}(P-P_0)<c_{5}^{-1}\max\{1,\hFplus(A/K)\},$$ (la hauteur \'etant le carr\'e de la norme) donc par la conjecture \ref{LangSilV3} il existe une sous-vari\'et\'e ab\'elienne stricte $B$ avec $\deg_L(B)\leq c_{5} \deg_L(A)$ et un entier naturel $1\leq N \leq c_{5}$ tel que $[N](P-P_0)\in B$. Tout point $P\in{C(K)}$ tel que $[N]P$ appartient \`a $B+[N]P_0$ est en fait un point vivant sur $[N]^{-1}(B+[N]P_0)\cap C$. On a bien s\^ur $\dim([N]^{-1}(B+[N]P_0)\cap C)\leq1$. Comme la courbe $C$ est plong\'ee dans sa jacobienne $\Jac(C)$ via $j$, la somme de $g$ copies $j(C)+\ldots+j(C)$ couvre tout $\Jac(C)$. Cette courbe ne peut donc pas, d\`es que $g\geq 2$, \^etre contenue dans un translat\'e de sous-vari\'et\'e ab\'elienne stricte. Donc $\dim([N]^{-1}(B+[N]P_0)\cap C) =0$ et on majore le degr\'e de l'intersection par le th\'eor\`eme de B\'ezout : $$\deg_{L}([N]^{-1}(B+[N]P_0)\cap C)\leq \deg_L([N]^{-1}(B+[N]P_0))\deg_L(C).$$ Par la proposition 2.3 page 9 de \cite{G} on obtient $$\deg_L([N]^{-1}(B+[N]P_0))=N^{2\codim(B)}\deg_L(B).$$ On d\'eduit donc $$\deg_{L}([N]^{-1}(B+[N]P_0)\cap C)\leq N^{2g}\deg_L(B)\deg_L(C).$$

On obtient alors la majoration suivante, en vertu de $\deg_L(B)\leq c_{5}\deg_L(A)$, de $N\leq c_{5}$ et de $\deg_L(C)=16g$ :
\begin{equation}\label{cardinal de check S1}
\#S_1\leq   \#I\;c_{5}^{2g+1}\deg_L(A)\; 16g
\end{equation}

\noindent et donc avec $\deg_L(A)=16^g g!$
\begin{equation}\label{cardinal de S1}
\#S_1\leq 16^{g+1}g!g\, c_{5}^{2g+1} \#I.
\end{equation}

On va \`a pr\'esent majorer le terme $\#I$ dans le membre de droite de (\ref{cardinal de S1}). En reprenant (\ref{cardinal}), le majorant de $c_{7}^2$ donn\'e dans le lemme \ref{c5} montre que 
\begin{equation}\label{split}
\#I\leq \left(1+\max\left\{2,\frac{4\cdot2^{98+5g}g^{11}\sqrt{\max\{1, h_{\Theta}(A)\}}}{\sqrt{c_{5}^{-1} \max\{1, \hFplus(A/K)\}}}\right\}\right)^{m_K}.
\end{equation}

Pour comparer la hauteur th\^eta et la hauteur de Faltings explicitement, on invoque \cite{D, T} afin d'obtenir 
\begin{equation}\label{comparaison}
\hFplus(A/K)\geq \hFplus(A/\overline{\mathbb{Q}})\geq 2h_{\Theta}(A)-\frac{g}{2}\log \max\{1,h_{\Theta}(A)\}-2M(g),
\end{equation}
o\`u $\displaystyle{M(g)=g4^{2g}}$ convient (voir le Corollaire \ref{hfaltheta} dans l'appendice du pr\'esent article).
On divise en deux cas.
\\

\underline{Premier cas :} si $h_{\Theta}(A)> 2M(g)+\frac{g}{2}\log\max\{1,h_{\Theta}(A)\},$ alors $$2h_{\Theta}(A)-2M(g)-\frac{g}{2}\log\max\{1,h_{\Theta}(A)\}>h_{\Theta}(A),$$ donc $\hFplus(A/K)> h_{\Theta}(A)$ et par (\ref{split}) on obtient (on rappelle que $c_5\geq 1$ puisque $c_5$ majore l'entier $N\geq1$) $$\#I \leq \left(1+\frac{2^{100+5g}g^{11}}{\sqrt{c_{5}^{-1}}}\right)^{m_K}.$$

\underline{Second cas :} si $h_{\Theta}(A)\leq 2M(g)+\frac{g}{2}\log\max\{1,h_{\Theta}(A)\},$ alors $h_{\Theta}(A)\leq 4M(g)$ et on utilise cela pour majorer le num\'erateur dans (\ref{split}) (tout en minorant le terme au d\'enominateur par $\sqrt{c_{5}^{-1}}$) et ainsi trouver
$$\#I\leq \left(1+2^{100+5g}g^{11}\sqrt{4c_{5}M(g)}\right)^{m_K}.$$
Dans les deux cas, en utilisant $M(g)=g4^{2g}$ (voir le Corollaire \ref{hfaltheta} dans l'appendice du pr\'esent article), on a montr\'e la majoration suivante (o\`u on utilise $c_{5}^{-1}\leq1$) :
\begin{equation}\label{cardinal de I}
\#I\leq \left(\sqrt{c_{5}}2^{102+7g}g^{12}\right)^{m_K}.
\end{equation}

On peut donc revenir \`a (\ref{cardinal de S1}) pour obtenir :

\begin{equation}\label{cardinal de S1 final}
\#S_1\leq 16^{g+1}g!g\, c_{5}^{2g+1} \Big(\sqrt{c_{5}} 2^{102+7g}g^{12}\Big)^{m_K}
\end{equation}

Combinant $(\ref{S2})$ et $(\ref{cardinal de S1 final})$, une derni\`ere majoration donne le r\'esultat.
\end{proof}

\section{Variante de Lang-Silverman}\label{varianteLS}

On regarde dans ce paragraphe une autre variante de Lang-Silverman, issue de \cite{V}.

\begin{conjecture}\label{LangSilV2}
Soit $g\geq 1$ un entier. Pour tout corps de nombres $K$, il existe trois nombres $c_{8}=c_{8}(g,K)>0$, $c_{9}=c_{9}(g,K)>0$ et $c_{10}=c_{10}(g,K)>0$  tels que pour toute vari\'et\'e ab\'elienne $A/K$ de dimension $g$ et pour tout fibr\'e en droites ample et sym\'etrique $L$ sur $A$, pour tout point $P\in{A(K)}$, au moins l'une des deux assertions suivantes est v\'erifi\'ee :
\begin{itemize}
\item[$\bullet$] il existe une sous-vari\'et\'e ab\'elienne $B\subset A$, $B\neq A$, de degr\'e $\deg_{L}(B)\leq c_{9}\deg_{L}(A)$ et telle que le point $P$ est d'ordre born\'e par $c_{10}$ modulo $B$,
\item[$\bullet$] on a $\mathrm{End}(A)\!\cdot\! P$ est Zariski dense et
\[
\widehat{h}_{A,L}(P) \geq c_{8}\, \max\Big\{\hFplus(A/K),\,1\Big\}\;,
\]
o\`u $\widehat{h}_{A,L}(.)$ est la hauteur de {N\'eron-Tate} associ\'ee \`a $L$ et $\hFplus(A/K)$ est la hauteur de Faltings de la vari\'et\'e $A/K$.
\end{itemize}
\end{conjecture}

C'est un \'enonc\'e qui dit sensiblement la m\^eme chose que la conjecture \ref{LangSilV3}: les points rationnels dont la hauteur est trop petite proviennent en fait d'un translat\'e d'une sous-vari\'et\'e ab\'elienne stricte. Cette version est l\`a encore motiv\'ee par le th\'eor\`eme 1.4 page 511 de \cite{F} et les th\'eor\`emes 1.8 et 1.13 de \cite{U}. Des exemples de familles de jacobiennes de courbes v\'erifiant cette conjecture sont donn\'es dans \cite{P} et \cite{U} : ces exemples sont aussi valables pour cette forme forte, puisque les jacobiennes associ\'ees sont simples.

Si on parcourt le raisonnement expos\'e plus haut, la m\^eme preuve vaut (en faisant de petites retouches ais\'ees, que nous faisons plus bas), le d\'ecompte est un peu moins pratique car il est \`a pr\'esent n\'ecessaire de savoir estimer le nombre de sous-vari\'et\'es ab\'eliennes  de degr\'e born\'e pour tenir compte de l'alternative, alors que dans la conjecture \ref{LangSilV3} n'y a qu'une seule vari\'et\'e obstructrice $B$.

\begin{remark}
Notons une diff\'erence entre la conjecture \ref{LangSilV3} et la conjecture \ref{LangSilV2}. Prenons le cas d'un point $P=(P_1,0)\in{A=A_1\times A_2}$ avec la polarisation $L$ produit de $L_1$ sur $A_1$ et $L_2$ sur $A_2$, avec $A_2$ ind\'ependante de $A_1$ (la hauteur de $A_2$ de d\'epend pas de la hauteur de $A_1$, la vari\'et\'e $A_2$ n'est donc pas dans la classe d'isog\'enie de $A_1$, ni isog\`ene \`a une puissance de $A_1$, etc). On suppose que $P_1$ n'est pas un point de torsion. On choisit $A_2$ telle que $\hFplus(A_2/K)>1+\widehat{h}_{A,L}(P)$, ce qui est possible puisque $\widehat{h}_{A,L}(P)=\widehat{h}_{A_1,L_1}(P_1)$ est ind\'ependant de $A_2$. Donc $P$ doit \^etre dans le premier cas. Regardons \`a pr\'esent l'ensemble des it\'er\'es de $P$, la vari\'et\'e $A_2$ \'etant fix\'ee. A partir d'un certain rang, la suite de r\'eels positifs $\widehat{h}_{A,L}([k]P)$ d\'epasse $\max\{1, \hFplus(A_2/K)\}$, et donc $[k]P$ passe dans le deuxi\`eme cas. C'est un ph\'enom\`ene qui ne se produit pas dans la conjecture \ref{LangSilV2} o\`u un point exclu du second cas voit tous ses it\'er\'es exclus \textit{de facto}. 
\end{remark}

\begin{remark}
La conjecture classique (la forme ``faible'') propose une minoration de la hauteur canonique d'un point rationnel par la hauteur de Faltings de la vari\'et\'e en imposant des hypoth\`eses n\'ecessaires sur le point, et ne dit rien sur les points ne satisfaisant pas ces hypoth\`eses, qui se retrouvent exclus. Un geste d'ouverture envers ces points a \'et\'e fait et la conjecture a ainsi \'et\'e fortifi\'ee pour la premi\`ere fois dans \cite{T}, l'inspiration venant du th\'eor\`eme 1.4 page 511 de \cite{F}, puis la forme forte a \'et\'e corrig\'ee dans \cite{V} pour donner la conjecture \ref{LangSilV2} et dans \cite{Y} pour donner la conjecture \ref{LangSilV3}. 
\end{remark}

Voyons \`a pr\'esent ce qui change dans le d\'ecompte si on fait usage de cette conjecture \ref{LangSilV2}.

\begin{theorem}\label{thm final 2}
Supposons vraie la conjecture $\ref{LangSilV2}$. Soit $K$ un corps de nombres et soit $g\geq 2$ un entier. Il existe une constante $c_{11}=c_{11}(g,K)>0$ telle que pour toute courbe $C$ de genre $g$ d\'efinie sur $K$, 
$$\#C(K)\leq c_{11}^{m_K+1},$$
o\`u $m_K$ d\'esigne le rang de Mordell-Weil de la jacobienne de $C$ sur $K$. On peut choisir
$$c_{11}=(2g)^{42g^4}\max\Big\{c_{9}^{(2g+1)4g^2+1}, c_{10}^{2g}, \frac{1}{\sqrt{c_{8}}}\Big\},$$ 
avec $c_{8}=c_{8}(g,K)>0$, $c_{9}=c_{9}(g,K)\geq 1$ et $c_{10}=c_{10}(g,K)>0$ donn\'ees dans la conjecture \ref{LangSilV2}.
\end{theorem}

\begin{proof}
Soit $\mathcal{B}$ l'ensemble des sous-vari\'et\'es ab\'eliennes de $A$ de degr\'e born\'e par $c_{9}\deg_L(A)$. C'est un ensemble fini, une borne sur son cardinal est donn\'ee un peu plus bas. On reprend la majoration \`a l'\'etape (\ref{cardinal de check S1}):
\begin{equation}
\#S_1\leq \deg_L(C) N^{2g}\sum_{P_0\in I}\,\sum_{B\in{\mathcal{B}}}\, \deg_L(B) ,
\end{equation}

\noindent donc comme pour tout $B\in{\mathcal{B}}$ on a $\deg_L(B)\leq c_{9}\deg_L(A)$, et comme $\deg_L(C)=16g$ et $N\leq c_{10}$  il vient
\begin{equation}\label{cardinal de check check S1}
\#S_1\leq 16g \;c_{10}^{2g}\; \#I\;\#\mathcal{B}\; c_{9}\deg_L(A) 
\end{equation}

\noindent et donc avec $\deg_L(A)=16^g g!$
\begin{equation}\label{cardinal S1}
\#S_1\leq 16^{g+1}g!g \;c_{10}^{2g}\; c_9\; \#I \;\#\mathcal{B}.
\end{equation}

L'\'equivalent de (\ref{cardinal de I}) est ici 
\begin{equation}\label{cardinal de hihi}
\#I\leq \left(\frac{1}{\sqrt{c_{8}}}2^{102+6g}g^{12}\right)^{m_K}.
\end{equation}

Majorons les autres termes du membre de droite de (\ref{cardinal S1}). Par la proposition 4.1 page 529 de \cite{W} on obtient la borne $\#\mathcal{B}\leq (4g!(2g)! (2^g c_{9}\deg_L(A))^{2g+1})^{4g^2}$, le majorant donn\'e ici \'etant l\'eg\`erement moins fin que celui de R\'emond, mais plus facile \`a manipuler. On a encore une fois $\deg_L(A)=16^g g!$, ce qui donne en majorant brutalement $\#\mathcal{B}\leq 2^{40g^4}g^{24g^4}c_9^{(2g+1)4g^2}$ et par (\ref{cardinal S1}) et (\ref{cardinal de hihi}) on obtient

\begin{equation}\label{cardinal de S1 final final}
\#S_1\leq  2^{41g^4} g^{25g^4} c_{10}^{2g}c_9^{(2g+1)4g^2+1} \Big(\frac{1}{\sqrt{c_8}} 2^{102+6g}g^{12}\Big)^{m_K}
\end{equation}

Combinant $(\ref{S2})$ et $(\ref{cardinal de S1 final final})$ et $g\geq 2$, un calcul facile donne le r\'esultat.
\end{proof}

On conclut avec une remarque suppl\'ementaire, donnant la mesure de la difficult\'e de la conjecture \ref{LangSilV2} (et la m\^eme remarque vaut pour la conjecture \ref{LangSilV3}).

\begin{remark}\label{borne torsion}
Supposons vraie la conjecture \ref{LangSilV2}. Alors pour tout $g\geq 1$, pour tout corps de nombres $K$, il existe une quantit\'e $c_{12}(g,K)>0$ telle que pour toute vari\'et\'e ab\'elienne $A$ de dimension $g$ et d\'efinie sur $K$, on a la majoration $\#A(K)_{tors} \leq c_{12}(g,K)$, o\`u la quantit\'e $c_{12}(g,K)$ est exprimable explicitement en termes des quantit\'es apparaissant dans la conjecture \ref{LangSilV2}.
\end{remark}

\section{Appendice: Corrigendum \`a ``Minorations des hauteurs normalis\'ees des sous-vari\'et\'es de vari\'et\'es ab\'eliennes II''}\label{DavPhi}

\begin{center}
par\\
\vspace{0.5cm}

{\Large\it Sinnou {\sc David} \&\ Patrice {\sc Philippon}
}
\end{center}
\vspace{0.5cm}

Une erreur s'est gliss\'ee dans la preuve du lemme 6.7 (lemme matriciel) et du corollaire 6.9 de l'appendice de~\cite{I}. Nous pr\'esentons ici une version corrig\'ee de ces \'enonc\'es, les constantes num\'eriques y sont chang\'ees. Nous indiquons d'embl\'ee que ces changements n'alt\`erent pas les autres \'enonc\'es de~\cite{I}. En effet, le lemme 6.7 n'est invoqu\'e que pour \'etablir le lemme suivant 6.8 et le corollaire 6.9 et c'est au lemme 6.8 qu'il est fait appel dans le texte. Or, la preuve du lemme 6.8 utilise en r\'ealit\'e une version faible de l'in\'egalit\'e du lemme 6.7, dont la v\'eracit\'e est confirm\'ee et que nous avons d\'etach\'ee dans l'\'enonc\'e~\ref{matrix} de la version corrig\'ee ci-dessous.

Nous remercions Pascal Autissier et Fabien Pazuki de nous avoir signal\'e cette erreur. Le lecteur pourra se reporter \`a~\cite{B} pour une version presque optimale du lemme matriciel; en effet le seul terme susceptible d'\^etre am\'elior\'e est le terme additif ne d\'ependant que de $g$. Pour une version accessible de la comparaison entre la hauteur th\^eta et la hauteur diff\'erentielle, le lecteur pourra consulter~\cite{T}. 

Nous reprenons ici les m\^emes conventions et notations que dans le texte initial. En particulier, pour tout \'el\'ement $\tau$ de l'espace de {\sc Siegel} ${\mathfrak S}_g$ nous notons ${k}_{\tau}$ le corps de d\'efinition de l'origine de $A_\tau$, c'est-\`a-dire du point de coordonn\'ees projectives $\left(\theta_p(\tau,0)\right)_{p\in{\mathcal{Z}}_2^2}$. On note de plus $F_g$ le domaine fondamental de Siegel tel que dans le texte initial. Rappelons enfin que, selon la notation~3.2 de~\cite{I}, $h(A)$ d\'esigne la hauteur projective de l'origine de $A$ dans le plongement induit par la puissance $16$-i\`eme de la polarisation principale implicite dans le choix de $\tau$.

\begin{lemma}\label{matrix}
Soit $\tau$ un \'el\'ement de ${\mathfrak S}_g$ tel que ${k}_{\tau}$ soit contenu dans un corps de nombres ${ k}$, de degr\'e $d$ sur $\mathbb{Q}$. Soit $\sigma$ un plongement complexe de ${ k}$ et $\tau(\sigma)$ un \'el\'ement de ${\mathfrak S}_g$ satisfaisant $A^{\sigma}_{\tau}=A_{\tau(\sigma)}$. Soit enfin $\tau'(\sigma)$ un repr\'esentant de $\tau(\sigma)$ dans $F_g$ et $y'(\sigma)$ la partie imaginaire de $\tau'(\sigma)$. Alors pour tout plongement $\sigma$ on a~:
$$\frac{\pi}{8}\Vert y'(\sigma)\Vert\leq
dh(A)+2g^2\log(4g)\enspace.$$
De plus,
$$\frac{\pi}{8d}\sum_{\sigma}\Vert y'(\sigma)\Vert\leq
4^gh(A)+2g^2\log(4g)\enspace.$$
\end{lemma}
\begin{proof} Quitte \`a faire une extension de ${ k}$, on peut supposer que $\tau'(\sigma)=\tau(\sigma)$; c'est ce que nous ferons. Soit donc $\sigma$ un plongement complexe de ${ k}$, choisissons des indices $p(\sigma),q(\sigma)\in{\mathcal{Z}}_{2}^2$ tels que
$\left|\theta_{p(\sigma)}(\tau({\sigma}),0)\right|$ soit minimal, non nul, tandis que $\left|\theta_{q(\sigma)}(\tau({\sigma}),0)\right|$ est maximal. En tenant compte
des lemmes~6.5 et~6.6 de~\cite{I} on a
\begin{equation}\label{estimations}
\left|\theta_{q(\sigma)}(\tau({\sigma}),0)\right|\geq1\quad\mbox{\rm et}\quad 0\neq \vert\theta_{p(\sigma)}(\tau({\sigma}),0)\vert\leq (4g)^{2g^2}\exp(-\frac{\pi}{8}\Vert y'({\sigma})\Vert)
\enspace.\end{equation}
Maintenant, pour $p\in{\mathcal{Z}}_2^2$ nous notons $I(p)$ l'ensemble des plongements tels que $p(\sigma)=p$.
La d\'efinition de la hauteur de {\sc Weil} nous assure alors que~:
\begin{equation}\label{premmaj}
\begin{array}{lcl}
\displaystyle
\frac{1}{d}\sum_{\sigma\in I(p)}\log\left\vert\frac{\theta_{q(\sigma)} (\tau({\sigma}), 0)}{\theta_{p}(\tau({\sigma}),0)}\right\vert 
& = & \displaystyle \frac{1}{d}\sum_{\sigma\in I(p)}\kern3pt\log\max_{q\in {\mathcal{Z}}_2^2}\left\{\left\vert\frac{\theta_{q}(\tau({\sigma}),
0)}{\theta_{p}(\tau({\sigma}),0)}\right\vert\right\}\\ 
& \leq & \displaystyle \frac{1}{d}\sum_{\sigma}\kern3pt\log\max_{q\in {\mathcal{Z}}_2^2}\left\{\left\vert\frac{\theta_{q}(\tau({\sigma}),
0)}{\theta_{p}(\tau({\sigma}),0)}\right\vert\right\}\\ 
& \leq & \rule{0mm}{8mm}\displaystyle h_{W,\kern1pt{\mathcal{L}}^{\otimes 16}}(A_{\tau}) = h(A)\enspace,
\end{array}
\end{equation}
car $\frac{\theta_{q}(\tau({\sigma}),0)}{\theta_{p}(\tau({\sigma}),0)} = \sigma\left(\frac{\theta_{q}(\tau,0)}{\theta_{p}(\tau,0)}\right)$ pour tous les plongements $\sigma$. On d\'eduit en particulier de~(\ref{estimations}) et~(\ref{premmaj}) pour tout plongement $\sigma$~:
$$\frac{\pi}{8}\Vert y'(\sigma)\Vert - 2g^2\log(4g) \leq \log\left\vert\frac{\theta_{q(\sigma)} (\tau({\sigma}), 0)}{\theta_{p(\sigma)}(\tau({\sigma}),0)}\right\vert \leq dh(A)
\enspace,$$
ce qui d\'emontre la premi\`ere in\'egalit\'e du lemme.

Par ailleurs, il suit de~(\ref{estimations})~:
$$\frac{1}{d}\sum_{\sigma}\log\left|\frac{\theta_{q(\sigma)}(\tau({\sigma}),
0)}{\theta_{p(\sigma)}(\tau({\sigma}),0)}\right| \geq \frac{\pi}{8d}\sum_{\sigma}\kern1pt\Vert y'({\sigma})\Vert-2g^2\log(4g)
\enspace.$$
Ainsi, le principe des tiroirs combin\'e \`a cette derni\`ere in\'egalit\'e montre qu'il existe un \'el\'ement $p\in{\mathcal{Z}}_2^2$ tel que
\begin{equation}\label{deuxmaj}
\begin{array}{lcl}
\displaystyle
\frac{\Card({\mathcal{Z}}_2^2)}{d}.\sum_{\sigma\in I(p)}\log\left|\frac{\theta_{q(\sigma)}(\tau({\sigma}),0)}{\theta_{p(\sigma)}(\tau({\sigma}),0)}\right| 
&\geq &\displaystyle\frac{1}{d}\sum_{\sigma}\log\left|\frac{\theta_{q(\sigma)}(\tau({\sigma}),
0)}{\theta_{p(\sigma)}(\tau({\sigma}),0)}\right|\\
&\geq &\displaystyle\frac{\pi}{8d}\sum_{\sigma}\Vert y'(\sigma)\Vert - 2g^2\log(4g)
\enspace.
\end{array}
\end{equation}
Reportant~(\ref{deuxmaj}) dans~(\ref{premmaj}) et observant $\Card({\mathcal{Z}}_2^2)=4^g$, on aboutit \`a la seconde in\'egalit\'e du lemme~\ref{matrix} qui est donc compl\`etement d\'emontr\'e.
\end{proof}

On notera que si l'on tient compte du th\'eor\`eme~1.1 de \cite{D}, {\em voir}~\cite{T}, th\'eor\`eme~1.1, on d\'eduit imm\'e\-diatement du lemme~\ref{matrix} et du lemme~6.4(i) de ~\cite{I} le corollaire suivant~:
\begin{corollary}\label{hfaltheta}
Soit $A$ une vari\'et\'e abelienne de dimension $g$, principalement polaris\'ee, d\'efinie sur un corps de nombres ${ k}$ sur lequel elle admet une r\'eduction semi-stable, et soit $h_{F}(A)$ sa hauteur de Faltings. Alors~:
$$\left\vert h(A)-\frac{1}{2} h_F(A)\right\vert\leq \frac{g}{4}\log\left(\max\left\{1;h(A)\right\}\right)+g4^{2g}\enspace.$$
\end{corollary}

\bibliographystyle{amsalpha}

\end{document}